\documentclass[11pt]{amsart}
\usepackage{amssymb}

\newtheorem{theorem}{Theorem}[section]
\newtheorem{corollary}[theorem]{Corollary}
\newtheorem{lemma}[theorem]{Lemma}
\newtheorem{proposition}[theorem]{Proposition}
\theoremstyle{definition}
\newtheorem{definition}[theorem]{Definition}

\newtheorem{question}[theorem]{Question}
\theoremstyle{remark}

\newtheorem{example}[theorem]{Example}

\numberwithin{equation}{section}

\begin{document}

\title{A Trichotomy for a Class of\\ Equivalence Relations}
\author{Longyun Ding}
\address{School of Mathematical Sciences and LPMC, Nankai University, Tianjin, 300071, P.R.China}
\email{dinglongyun@gmail.com}
\thanks{Research partially supported by the National Natural Science Foundation of China (Grant No. 10701044).}


\date{\today}

\begin{abstract}
Let $X_n,\,n\in\Bbb N$ be a sequence of non-empty sets,
$\psi_n:X_n^2\to\Bbb R^+$. We consider the relation
$E((X_n,\psi_n)_{n\in\Bbb N})$ on $\prod_{n\in\Bbb N}X_n$ by
$(x,y)\in E((X_n,\psi_n)_{n\in\Bbb N})\Leftrightarrow\sum_{n\in\Bbb
N}\psi_n(x(n),y(n))<+\infty$. If $E((X_n,\psi_n)_{n\in\Bbb N})$ is a
Borel equivalence relation, we show a trichotomy that either $\Bbb
R^\Bbb N/\ell_1\le_B E$, $E_1\le_B E$, or $E\le_B E_0$.

We also prove that, for a rather general case,
$E((X_n,\psi_n)_{n\in\Bbb N})$ is an equivalence relation iff it is
an $\ell_p$-like equivalence relation.
\end{abstract}
\maketitle

\section{Introduction}

A topological space is called a {\it Polish space} if it is
separable and completely metrizable. Let $X,Y$ be Polish spaces and
$E,F$ equivalence relations on $X,Y$ respectively. A {\it Borel
reduction} of $E$ to $F$  is a Borel function $\theta:X\to Y$ such
that $(x,y)\in E$ iff $(\theta(x),\theta(y))\in F$, for all $x,y\in
X$. We say that $E$ is {\it Borel reducible} to $F$, denoted $E\le_B
F$, if there is a Borel reduction of $E$ to $F$. If $E\le_B F$ and
$F\le_B E$, we say that $E$ and $F$ are {\it Borel bireducible} and
denote $E\sim_B F$. We refer to \cite{gao} and \cite{kanovei} for
background on Borel reducibility.

There are several famous dichotomy theorems on Borel reducibility.
The first one is the Silver's dichotomy theorem \cite{silver}.

\begin{theorem}[Silver]
Let $E$ be a ${\bf\Pi}^1_1$ equivalence relation. Then $E$ has
either at most countably many or perfectly many equivalence classes,
i.e. $E\le_B{\rm id}(\Bbb N)\mbox{ or }{\rm id}(\Bbb R)\le_B E$.
\end{theorem}

There are three dichotomy theorems concerning $E_0$. Before
introducing these theorems, we recall definitions of equivalence
relations $E_0,E_1,E_0^\omega$.
\begin{enumerate}
\item[(a)] For $x,y\in 2^\Bbb N$, $(x,y)\in E_0\Leftrightarrow\exists m\forall
n\ge m(x(n)=y(n)).$
\item[(b)] For $x,y\in 2^{\Bbb N\times\Bbb N}$, $(x,y)\in
E_1\Leftrightarrow\exists m\forall n\ge m\forall k(x(n,k)=y(n,k)).$
\item[(c)] For $x,y\in 2^{\Bbb N\times\Bbb N}$, $(x,y)\in
E_0^\omega\Leftrightarrow\forall k\exists m\forall n\ge
m(x(n,k)=y(n,k)).$
\end{enumerate}

\begin{theorem}
Let $E$ be a Borel equivalence relation. Then
\begin{enumerate}
\item[(a)] {\rm (Harrington-Kechris-Louveau \cite{HKL})} either $E\le_B{\rm
id}(\Bbb R)$ or $E_0\le_B E$;
\item[(b)] {\rm (Kechris-Louveau \cite{KL})} if $E\le_B E_1$, then
$E\le_B E_0$ or $E\sim_B E_1$;
\item[(c)] {\rm (Hjorth-Kechris \cite{HK})} if $E\le_B E_0^\omega$,
then $E\le_B E_0$ or $E\sim_B E_0^\omega$.
\end{enumerate}
\end{theorem}

Another class of interesting Borel equivalence relations come from
classical Banach sequence spaces. Let $p\ge 1$. For $x,y\in\Bbb
R^\Bbb N$, $(x,y)\in\Bbb R^\Bbb N/\ell_p\Leftrightarrow
x-y\in\ell_p$. It was shown by G. Hjorth \cite{hjorth} that every
Borel equivalence relation $E\le_B\Bbb R^\Bbb N/\ell_1$ is either
essentially countable or satisfies $E\sim_B\Bbb R^\Bbb N/\ell_1$.
Kanovei asked whether the position of $\Bbb R^\Bbb N/\ell_1$ in the
$\le_B$-structure is similar with $E_1$ and $E_0^\omega$ (see
\cite{kanovei}, Question 5.7.5).

\begin{question}[Kanovei]
Does every Borel equivalence relation $E\le\Bbb R^\Bbb N/\ell_1$
satisfy either $E\le_B E_0$ or $E\sim_B\Bbb R^\Bbb N/\ell_1$?
\end{question}

Two kinds of $\ell_p$-like equivalence relations were introduced by
T. M\' atrai \cite{matrai} and the author \cite{ding}. (1) Let
$f:[0,1]\to\Bbb R^+$. For $x,y\in[0,1]^\Bbb N$, $(x,y)\in
E_f\Leftrightarrow\sum_{n\in\Bbb N}f(|x(n)-y(n)|)<+\infty$. (2) Let
$(X_n,d_n),\,n\in\Bbb N$ be a sequence of metric spaces, $p\ge 1$.
For $(x,y)\in\prod_{n\in\Bbb N}X_n$, $(x,y)\in E((X_n,d_n)_{n\in\Bbb
N};p)\Leftrightarrow\sum_{n\in\Bbb N}d_n(x(n),y(n))^p<+\infty$.

In this paper, we introduce a notion surpassing both (1) and (2).
Let $X_n,\,n\in\Bbb N$ be a sequence of non-empty sets,
$\psi_n:X_n^2\to\Bbb R^+$. For $x,y\in\prod_{n\in\Bbb N}X_n$,
$(x,y)\in E((X_n,\psi_n)_{n\in\Bbb N})\Leftrightarrow\sum_{n\in\Bbb
N}\psi_n(x(n),y(n))<+\infty$. Though we did not find a natural
necessary and sufficient condition that $E((X_n,\psi_n)_{n\in\Bbb
N})$ be an equivalence relation, we establish the following
trichotomy.

\begin{theorem}
If $E=E((X_n,\psi_n)_{n\in\Bbb N})$ is a Borel equivalence relation,
then either $\Bbb R^\Bbb N/\ell_1\le_B E$, $E_1\le_B E$, or $E\le_B
E_0$.
\end{theorem}

From this trichotomy, we can see that Kanovei's problem is valid
within equivalence relations of the form $E((X_n,\psi_n)_{n\in\Bbb
N})$.

It was shown by R. Dougherty and G. Hjorth \cite{DH} that, for
$p,q\ge 1$, $\Bbb R^\Bbb N/\ell_p\le_B\Bbb R^\Bbb N/\ell_q$ iff
$p\le q$. We complete a different picture for $0<p\le 1$ by showing
that $\Bbb R^\Bbb N/\ell_p\sim_B\Bbb R^\Bbb N/\ell_1$.

Via a process of metrization, we prove that, for a rather general
case, equivalence relations $E((X_n,\psi_n)_{n\in\Bbb N})$ coincide
with $\ell_p$-like equivalence relations $E((X_n,d_n)_{n\in\Bbb
N};p)$.

\section{A trichotomy for sum-like equivalence relations}

We denote the set of all non-negative real numbers by $\Bbb R^+$.

\begin{definition}
Let $X_n,\,n\in\Bbb N$ be a sequence of non-empty sets,
$\psi_n:X_n^2\to\Bbb R^+$. We define a relation
$E((X_n,\psi_n)_{n\in\Bbb N})$ on $\prod_{n\in\Bbb N}X_n$ by
$$(x,y)\in E((X_n,\psi_n)_{n\in\Bbb N})\iff\sum_{n\in\Bbb
N}\psi_n(x(n),y(n))<+\infty$$ for $x,y\in\prod_{n\in\Bbb N}X_n$. If
$(X_n,\psi_n)=(X,\psi)$ for every $n\in\Bbb N$, we write
$E(X,\psi)=E((X_n,\psi_n)_{n\in\Bbb N})$ for the sake of brevity.
\end{definition}

It is hard to find a natural necessary and sufficient condition to
determine when $E((X_n,\psi_n)_{n\in\Bbb N})$ is an equivalence
relation. We need the following definition:

\begin{definition}
If $E((X_n,\psi_n)_{n\in\Bbb N})$ is an equivalence relation, we
call it a {\it sum-like equivalence relation}. Furthermore, if
$X_n,\,n\in\Bbb N$ is a sequence of Polish spaces and every $\psi_n$
is Borel function, it is called a {\it sum-like Borel equivalence
relation}.
\end{definition}

The following easy lemma is very useful for the study of sum-like
equivalence relations.

\begin{lemma} \label{ne}
Let $E=E((X_n,\psi_n)_{n\in\Bbb N})$ be a sum-like equivalence
relation. For $x,y\in\prod_{n\in\Bbb N}X_n$, we have
$$(x,y)\in E\iff\sum_{x(n)\ne y(n)}\psi_n(x(n),y(n))<+\infty.$$
\end{lemma}

\begin{proof}
Since $E$ is an equivalence relation, $(x,x)\in E$. It follows that
$$\sum_{x(n)=y(n)}\psi_n(x(n),y(n))\le\sum_{n\in\Bbb
N}\psi_n(x(n),x(n))<+\infty.$$ Then the lemma follows.
\end{proof}

\begin{definition}
Let $(X_n,d_n),\,n\in\Bbb N$ be a sequence of pseudo-metric spaces.
For $p\ge 1$, an {\it $\ell_p$-like equivalence relation}
$E((X_n,d_n)_{n\in\Bbb N};p)$ is $E((X_n,\psi_n)_{n\in\Bbb N})$
where $\psi_n(u,v)=d_n(u,v)^p$ for $u,v\in X_n$. If
$(X_n,d_n)=(X,d)$ for every $n\in\Bbb N$, we write
$E(X,d;p)=E((X_n,d_n)_{n\in\Bbb N};p)$ for the sake of brevity.
\end{definition}

\begin{proposition}
Let $E((X_n,d_n)_{n\in\Bbb N};p)$ be an $\ell_p$-like equivalence
relation. Then there exists metric $d_n'$ on $X_n$ for each $n$ such
that $E((X_n,d_n')_{n\in\Bbb N};p)=E((X_n,d_n)_{n\in\Bbb N};p)$.
\end{proposition}

\begin{proof}
We can see that following $d_n'$'s meet the requirements.
$$d_n'(u,v)=\left\{\begin{array}{ll}0, & u=v,\cr 2^{-n}, &
u\ne v,d_n(u,v)\le 2^{-n},\cr d_n(u,v), &
d_n(u,v)>2^{-n},\end{array}\right.$$ for $u,v\in X_n$.
\end{proof}

\begin{proposition}
Let $E_i,\,i\in\Bbb N$ be a sequence of Borel equivalence relations.
Then for every $p\ge 1$ there is an $\ell_p$-like Borel equivalence
relation $E=E((X_n,d_n)_{n\in\Bbb N};p)$ such that $E_i\le_B E$ for
$n\in\Bbb N$.
\end{proposition}

\begin{proof}
For $i\in\Bbb N$, let $E_i$ be a Borel equivalence relation on a
Polish space $Y_i$. We define $\chi_i:Y_i^2\to\Bbb R^+$ by
$$\chi_i(u,v)=\left\{\begin{array}{ll}0, &(u,v)\in E_i,\cr 1,
&(u,v)\notin E_i.\end{array}\right.$$ Now fix a bijection
$\langle\cdot,\cdot\rangle:\Bbb N^2\to\Bbb N$. For every $i,j\in\Bbb
N$, if $n=\langle i,j\rangle$, we denote $X_n=Y_i$ and $d_n=\chi_i$.
It is easy to verify that $E=E((X_n,d_n)_{n\in\Bbb N};p)$ is an
$\ell_p$-like Borel equivalence relation. For $i\in\Bbb N$, define
$\theta_i:Y_i\to\prod_{n\in\Bbb N}X_n$ by
$$\theta_i(u)(\langle k,j\rangle)=\left\{\begin{array}{ll}u, &
k=i,\cr a_k, & k\ne i,\end{array}\right.$$ where $a_k\in Y_k$ is
independent of $u$. Clearly $\theta_i$ is a Borel reduction of $E_i$
to $E$.
\end{proof}

Now suppose that $E=E((X_n,\psi_n)_{n\in\Bbb N})$ is a sum-like
Borel equivalence relation. Its position in the $\le_B$-structure of
Borel equivalence relations is determined by the following
condition:

$(\ell 1)$ {\it $\forall c>0\exists x,y\in\prod_{n\in\Bbb N}X_n$
such that $\exists m\forall n>m(\psi_n(x(n),y(n))<c)$ and
$$\sum_{n\in\Bbb N}\psi_n(x(n),y(n))=+\infty.$$}

\begin{lemma}
If $(\ell 1)$ holds, then $\Bbb R^\Bbb N/\ell_1\le_B E$.
\end{lemma}

\begin{proof}
Firstly, we show a well known fact: $\Bbb R^\Bbb
N/\ell_1\le_B[0,1]^\Bbb N/\ell_1$. Fix a bijection
$\langle\cdot,\cdot\rangle':\Bbb N\times\Bbb Z\to\Bbb N$. For
$z\in\Bbb R^\Bbb N$, we define $\theta'(z)\in[0,1]^\Bbb N$ as
$$\theta'(z)(\langle m,k\rangle')=\left\{\begin{array}{ll}0,&
z(m)<k,\cr z(m)-k,& k\le z(m)<k+1,\cr 1,& k+1\le
z(m).\end{array}\right.$$ Then $\theta'$ witness that $\Bbb R^\Bbb
N/\ell_1\le_B[0,1]^\Bbb N/\ell_1$.

Secondly, we construct a reduction of $[0,1]^\Bbb N/\ell_1$ to $E$.

For $l\in\Bbb N$, by condition $(\ell 1)$, there exist
$x_l,y_l\in\prod_{n\in\Bbb N}X_n$ such that $\exists m\forall
n>m(\psi_n(x_l(n),y_l(n))<2^{-l})$ and $\sum_{n\in\Bbb
N}\psi_n(x_l(n),y_l(n))=+\infty$. Then we can select two sequences
of natural numbers $(i_l)_{l\in\Bbb N},(j_l)_{l\in\Bbb N}$
satisfying that
\begin{enumerate}
\item[(i)] $i_l<j_l<i_{l+1}$ for $l\in\Bbb N$;
\item[(ii)] $\psi_n(x_l(n),y_l(n))<2^{-l}$ for $i_l\le n\le j_l$;
\item[(iii)] $1\le\sum_{i_l\le n\le j_l}\psi_n(x_l(n),y_l(n))<1+2^{-l}$.
\end{enumerate}

Fix an element $a_n\in X_n$ for every $n\in\Bbb N$. For
$z\in[0,1]^\Bbb N$, we define $\vartheta(z)\in\prod_{n\in\Bbb N}X_n$
by
$$\vartheta(z)(n)=\left\{\begin{array}{ll}x_l(n),& i_l\le
n\le j_l,\sum_{i_l\le m\le n}\psi_{m}(x_l(m),y_l(m))\le z(l);\cr
y_l(n),& i_l\le n\le j_l,\sum_{i_l\le m\le
n}\psi_{m}(x_l(m),y_l(m))>z(l);\cr a_n,&
\mbox{otherwise}.\end{array}\right.$$

Note that, for $z,w\in[0,1]^\Bbb N$ and $l\in\Bbb N$, we have
$$|z(l)-w(l)|-2^{-l}<\hspace{-20pt}\sum_{\mbox{\scriptsize$\begin{array}{c}i_l\le n\le j_l\\
\vartheta(z)(n)\ne\vartheta(w)(n)\end{array}$}}\hspace{-25pt}\psi_n(\vartheta(z)(n),\vartheta(w)(n))<|z(l)-w(l)|+2^{-l}.$$
Therefore, by Lemma \ref{ne},
$$\begin{array}{ll}(\vartheta(z),\vartheta(w))\in
E&\iff\sum_{\vartheta(z)(n)\ne\vartheta(w)(n)}\psi_n(\vartheta(z)(n),\vartheta(w)(n))<+\infty\cr
&\iff\sum_{l\in\Bbb N}|z(l)-w(l)|<+\infty\cr &\iff
z-w\in\ell_1.\end{array}$$ It follows that $[0,1]^\Bbb N/\ell_1\le_B
E$.
\end{proof}

If $(\ell 1)$ fails,  then there exists $c>0$ such that
$$\exists m\forall n>m(\psi_n(x(n),y(n))<c)\Rightarrow\sum_{n\in\Bbb
N}\psi_n(x(n),y(n))<+\infty$$ for $x,y\in\prod_{n\in\Bbb N}X_n$. Now
we denote
$$F_n=\{(u,v)\in X_n^2:\psi_n(u,v)<c\}.$$

\begin{lemma} \label{fn}
If $(\ell 1)$ fails, let $c>0$ and $F_n,\,n\in\Bbb N$ be defined as
above. Then
\begin{enumerate}
\item[(1)] for $x,y\in\prod_{n\in\Bbb N}X_n$,
$$(x,y)\in E\iff\exists m\forall n>m((x(n),y(n))\in F_n);$$
\item[(2)] there exists $N_0$ such that, for $n>N_0$, $F_n$ is a Borel equivalence
relation on $X_n$.
\end{enumerate}
\end{lemma}

\begin{proof} Since $(\ell 1)$ fails, clause (1) is trivial.

(2) Firstly, assume for contradiction that, there exist a strictly
increasing sequence of natural numbers $(n_k)_{k\in\Bbb N}$ and
$u_k\in X_{n_k}$ such that $\psi_{n_k}(u_k,u_k)\ge c$. Select an
$x\in\prod_{n\in\Bbb N}X_n$ with $x(n_k)=u_k$ for $k\in\Bbb N$. Then
we have $\sum_{n\in\Bbb N}\psi_n(x(n),x(n))=+\infty$. This is
impossible, since $E$ is an equivalence relation. So there exists
$N_1$ such that, for $n>N_1,u\in X_n$, we have $\psi_n(u,u)< c$,
i.e. $(u,u)\in F_n$.

Secondly, assume for contradiction that, there exist a strictly
increasing sequence of natural numbers $(n_k)_{k\in\Bbb N}$ and
$u_k,v_k\in X_{n_k}$ such that
$\psi_{n_k}(u_k,v_k)<c,\psi_{n_k}(v_k,u_k)\ge c$. Select
$x,y\in\prod_{n\in\Bbb N}X_n$ such that $x(n_k)=u_k,y(n_k)=v_k$ for
$k\in\Bbb N$ and $x(n)=y(n)$ for other $n$. We have
$\psi_n(x(n),y(n))<c$ for $n>N_1$. Since $(\ell 1)$ fails,
$\sum_{n\in\Bbb N}\psi_n(x(n),y(n))<+\infty$, i.e. $(x,y)\in E$.
Clearly $\sum_{n\in\Bbb N}\psi_n(y(n),x(n))=+\infty$, so
$(y,x)\notin E$. A contradiction! Hence there exists $N_2$ such
that, for $n>N_2,u,v\in X_n$, we have
$\psi_n(u,v)<c\Rightarrow\psi_n(v,u)<c$, i.e. $(u,v)\in
F_n\Rightarrow(v,u)\in F_n$.

With a similar argument, we can prove that there exists $N_3$ such
that, for $n>N_3, u,v,r\in X_n$, we have $(u,v),(v,r)\in
F_n\Rightarrow(u,r)\in F_n$.

In summary, for $n>N_0=\max\{N_1,N_2,N_3\}$, $F_n$ is an equivalence
relation. Since $\psi_n$ is Borel, so is $F_n$.
\end{proof}

Recall that $E_0(\Bbb N)$ is an equivalence relation on $\Bbb N^\Bbb
N$ similar to $E_0$ on $2^\Bbb N$. For $x,y\in\Bbb N^\Bbb N$,
$(x,y)\in E_0(\Bbb N)\Leftrightarrow\exists m\forall
n>m(x(n)=y(n))$. It is well known that $E_0\sim_B E_0(\Bbb N)$ (see
Proposition 6.1.2 of \cite{gao}).

\begin{lemma}
For any equivalence relation $E$ on $\prod_{n\in\Bbb N}X_n$, if
there is a sequence $F_n\subseteq X_n^2,\,n\in\Bbb N$ such that (1)
and (2) of Lemma \ref{fn} hold, then either $E_1\le_B E$, $E\sim_B
E_0$, or $E$ is trivial, i.e. all elements in $\prod_{n\in\Bbb
N}X_n$ are equivalent.
\end{lemma}

\begin{proof}
For $n>N_0$, since $F_n$ is a Borel equivalence relation on $X_n$,
from the Silver dichotomy theorem \cite{silver}, either there are at
most countably many $F_n$-equivalence classes or there are perfectly
many $F_n$-equivalence classes.

{\it Case 1.} There exists a strictly increasing sequence of natural
numbers $(n_k)_{k\in\Bbb N}$ such that there are perfectly many
$F_{n_k}$-equivalence classes. Then there is a continuous embedding
$h_{k}:2^\Bbb N\to X_{n_k}$ for every $k\in\Bbb N$ such that
$(h_k(z),h_k(w))\in F_{n_k}$ iff $z=w$. Define $\theta:2^{\Bbb
N\times\Bbb N}\to\prod_{n\in\Bbb N}X_n$ by
$$\theta(x)(n)=\left\{\begin{array}{ll}h_k(x(k,\cdot)),& n=n_k,\cr
a_n,&\mbox{otherwise,}\end{array}\right.$$ where $a_n\in X_n$ is
independent of $x$. By Lemma \ref{fn}.(1), it is straightforward to
check that $\theta$ is a reduction of $E_1$ to $E$.

{\it Case 2.} There exists $N$ such that, for $n>N$, $F_n$ has only
one equivalence class. From Lemma \ref{fn}.(1), we see that $E$ is
trivial.

{\it Case 3.} If case 1 fails, then there exists $N'$ such that, for
$n>N'$, $F_n$ has at most countably many equivalence classes. So
$E\le_B E_0(\Bbb N)$. If case 2 fails, then there exists a strictly
increasing sequence of natural numbers $(n_k)_{k\in\Bbb N}$ such
that $F_{n_k}$ has more than one equivalence class. Thus $E_0\le_B
E$. Since $E_0\sim_B E_0(\Bbb N)$, we have $E\sim_B E_0$.
\end{proof}

Now we have already completed the proof of the following trichotomy.

\begin{theorem}
Let $E=E((X_n,\psi_n)_{n\in\Bbb N})$ be a sum-like Borel equivalence
relation. Then either $\Bbb R^\Bbb N/\ell_1\le_B E$, $E_1\le_B E$,
or $E\le_B E_0$.
\end{theorem}

\begin{corollary}
Let $E=E((X_n,\psi_n)_{n\in\Bbb N})$ be a sum-like Borel equivalence
relation. If $E\le_B\Bbb R^\Bbb N/\ell_1$, then either $E\le_B E_0$
or $E\sim_B\Bbb R^\Bbb N/\ell_1$.
\end{corollary}

\begin{proof}
It is well known that $E_1\not\le_B\Bbb R^\Bbb N/\ell_1$ (see
\cite{KL} Theorem 4.2), so if $E\le_B\Bbb R^\Bbb N/\ell_1$, then
$E_1\not\le_B E$. Hence the corollary follows.
\end{proof}

It was shown by R. Dougherty and G. Hjorth \cite{DH} that, for
$p,q\ge 1$,
$$\Bbb R^\Bbb N/\ell_p\le_B\Bbb R^\Bbb N/\ell_q\iff p\le q.$$
The following corollary shows that, for $p\le 1$, the situation is
different.

\begin{corollary}
For $0<p\le 1$, we have $\Bbb R^\Bbb N/\ell_p\sim_B\Bbb R^\Bbb
N/\ell_1$.
\end{corollary}

\begin{proof}
Note that $\Bbb R^\Bbb N/\ell_p=E(\Bbb R,\psi)$ where
$\psi(u,v)=|u-v|^p$ for $u,v\in\Bbb R$. It is easy to see that
$(\ell 1)$ holds for $E(\Bbb R,\psi)$, so $\Bbb R^\Bbb
N/\ell_1\le_B\Bbb R^\Bbb N/\ell_p$.

For the other direction, we claim that $\Bbb R^\Bbb
N/\ell_p\le_B\Bbb R^\Bbb N/\ell_q$ for $0<p<q\le 1$.

We sketch the proof for Theorem 1.1 of \cite{DH}, that $\Bbb R^\Bbb
N/\ell_p\le_B\Bbb R^\Bbb N/\ell_q$ for $1\le p<q$, and check that it
is also valid for $0<p<q\le 1$.

It will suffice to prove for the case $0<\frac{q}{2}<p<q\le 1$. We
denote $\rho=\frac{p}{q}$ and $r=4^{-\rho}$. Then
$\frac{1}{2}<\rho<1$ and $\frac{1}{4}<r<\frac{1}{2}$. In \cite{DH}
p. 1838, the authors constructed a continuous function $\bar
K_r:\Bbb R\to\Bbb R^2$, and proved that there are $m',M'>0$ such
that
$$m'|s-t|^\rho\le\|\bar K_r(s)-\bar K_r(t)\|_2\le M'|s-t|^\rho,$$
for $s,t\in[i-1,i+1],i\in\Bbb Z$. And $\|\bar K_r(s)-\bar
K_r(t)\|_2\ge 1$ if $s,t$ are not in the same interval
$[i-1,i+1],i\in\Bbb Z$.

Define a mapping $\theta:\Bbb R^\Bbb N\to\Bbb R^\Bbb N$ such that,
for $x\in\Bbb R^\Bbb N$ and $k\in\Bbb N$,
$$\bar K_r(x(k))=(\theta(x)(2k),\theta(x)(2k+1)).$$
For $w=(s,t)\in\Bbb R^2$, denote
$\|w\|_q=(|s|^q+|t|^q)^\frac{1}{q}$. Note that
$$\frac{1}{\sqrt{2}}\|w\|_2\le\|w\|_\infty\le\|w\|_q\le
2^\frac{1}{q}\|w\|_\infty\le 2^\frac{1}{q}\|w\|_2.$$ For $x,y\in\Bbb
R^\Bbb N$, we have
$$\begin{array}{ll}\theta(x)-\theta(y)\in\ell_q&\iff\sum_{k\in\Bbb N}\|\bar K_r(x(k))-\bar K_r(y(k))\|_q^q<+\infty\cr
&\iff\sum_{k\in\Bbb N}\|\bar K_r(x(k))-\bar
K_r(y(k))\|_2^q<+\infty\cr &\iff\sum_{k\in\Bbb
N}|x(k)-y(k)|^p<+\infty\cr &\iff x-y\in\ell_p.\end{array}$$ Thus,
$\theta$ is a reduction of $\Bbb R^\Bbb N/\ell_p$ to $\Bbb R^\Bbb
N/\ell_q$.
\end{proof}

\section{Metrization}

In this section, we show that a sum-like equivalence relation
$E((X_n,\psi_n)_{n\in\Bbb N})$ coincides with an $\ell_p$-like
equivalence relation if the following conditions hold.

$(m1)$ {\it Denote $X=\bigcup_{n\in\Bbb N}X_n$. There is a unique
function $\psi:X^2\to\Bbb R^+$ such that $\psi_n=\psi\upharpoonright
X_n^2$ and $\psi(u,v)=1$ if no $X_n$ contains both $u,v$.}

$(m2)$ {\it For any $u,v,r\in X$, if $u,v,r\in X_n$, then there
exists $m>n$ such that $u,v,r\in X_m$.}

Let $\overline\psi_n=\min\{\psi_n,1\}$. We can see that
$E((X_n,\overline\psi_n)_{n\in\Bbb N})=E((X_n,\psi_n)_{n\in\Bbb
N})$. Therefore, we may assume $\psi_n\le 1$ if needed.

\begin{lemma} \label{C}
Let $(X_n,\psi_n)_{n\in\Bbb N}$ satisfy $(m1)$ and $(m2)$. If
$\psi_n\le 1$, then $E((X_n,\psi_n)_{n\in\Bbb N})$ is an equivalence
relation iff the following conditions hold:
\begin{enumerate}
\item[(i)] $\psi(u,u)=0$ for $u\in X$;
\item[(ii)] there is a $C\ge 1$ such that for $u,v,r\in X$,
$$\psi(v,u)\le C\psi(u,v);\quad\psi(u,r)\le
C(\psi(u,v)+\psi(v,r)).$$
\end{enumerate}
\end{lemma}

\begin{proof}
If conditions (i) and (ii) hold, it is trivial that
$E((X_n,\psi_n)_{n\in\Bbb N})$ is an equivalence relation. We only
need to prove the other direction.

Now assume that $E((X_n,\psi_n)_{n\in\Bbb N})$ is an equivalence
relation.

Firstly, for any $u\in X$, by $(m2)$, there is an infinite set
$I\subseteq\Bbb N$ such that $u\in X_n$ for $n\in I$. Let
$x\in\prod_{n\in\Bbb N}X_n$ with $x(n)=u$ for $n\in I$. Since
$(x,x)\in E((X_n,\psi_n)_{n\in\Bbb N})$, we have $\psi(u,u)=0$.

Secondly, for $u,v\in X$, if $\psi(u,v)=0$, we claim that
$\psi(v,u)=0$. By $(m2)$, there is an infinite set $I\subseteq\Bbb
N$ such that $u,v\in X_n$ for $n\in I$. Therefore, for any
$x,y\in\prod_{n\in\Bbb N}X_n$, if $x(n)=u,y(n)=v$ for $n\in I$ and
$x(n)=y(n)$ for $n\notin I$, we have $(x,y)\in
E((X_n,\psi_n)_{n\in\Bbb N})$. Hence $(y,x)\in
E((X_n,\psi_n)_{n\in\Bbb N})$. It follows that $\sum_{n\in
I}\psi(v,u)=\sum_{n\in I}\psi(y(n),x(n))<+\infty$. Thus
$\psi(v,u)=0$.

Now assume for contradiction that, for every $k\in\Bbb N$ there are
$u_k,v_k\in X$ such that $\psi(v_k,u_k)>2^k\psi(u_k,v_k)>0$. Since
$\psi_n\le 1$, $0<\psi(u_k,v_k)<2^{-k}$. By $(m2)$, there are
infinitely many $n$ such that $u_k,v_k\in X_n$.

Select a finite set $I_k\subseteq\Bbb N$ for every $k$ satisfying
that
\begin{enumerate}
\item[(i)] $u_k,v_k\in X_n$ for $n\in I_k$;
\item[(ii)] $2^{-k}\le|I_k|\psi(u_k,v_k)\le 2^{-(k-1)}$;
\item[(iii)] if $k_1<k_2$, then $\max I_{k_1}<\min I_{k_2}$.
\end{enumerate}

Now we define $x,y\in\prod_{n\in\Bbb N}X_n$ by
$$\left\{\begin{array}{ll}x(n)=u_k,y(n)=v_k,& n\in I_k,k\in\Bbb N,\cr
x(n)=y(n)=a_n,&\mbox{otherwise},\end{array}\right.$$ where $a_n\in
X_n$ is independent of $x$ and $y$. Then we have
$$\sum_{n\in\Bbb N}\psi(x(n),y(n))=\sum_{k\in\Bbb
N}|I_k|\psi(u_k,v_k)\le\sum_{k\in\Bbb N}2^{-(k-1)}<+\infty,$$ so
$(x,y)\in E((X_n,\psi_n)_{n\in\Bbb N})$. On the other hand, we have
$$\sum_{n\in\Bbb
N}\psi(y(n),x(n))=\sum_{k\in\Bbb N}|I_k|\psi(v_k,u_k)>\sum_{k\in\Bbb
N}2^k|I_k|\psi(u_k,v_k)\ge\sum_{k\in\Bbb N}1=+\infty,$$ so
$(y,x)\notin E((X_n,\psi_n)_{n\in\Bbb N})$. A Contradiction! We
complete the proof of that there is $C_1\ge 1$ such that
$\psi(u,v)\le C_1\psi(v,u)$ for $u,v\in X$.

With a similar argument, we can prove that there is $C_2\ge 1$ such
that $\psi(u,r)\le C_2(\psi(u,v)+\psi(v,r))$ for $u,v,r\in X$.
\end{proof}

\begin{lemma}
Let $E((X_n,\psi_n)_{n\in\Bbb N}),E((X_n,\varphi_n)_{n\in\Bbb N})$
be two sum-like equivalence relations, both satisfying $(m1)$ and
$(m2)$. If $\varphi_n\le 1$, then
$$\begin{array}{cl}& E((X_n,\psi_n)_{n\in\Bbb
N})\subseteq E((X_n,\varphi_n)_{n\in\Bbb N})\cr\iff& \exists A\ge
1\forall u,v\in X(\varphi(u,v)\le A\psi(u,v)).\end{array}$$
\end{lemma}

\begin{proof}
``$\Leftarrow$'' is trivial. ``$\Rightarrow$'' follows similarly as
the proof of Lemma \ref{C}.
\end{proof}

\begin{corollary}
Let $E((X_n,\psi_n)_{n\in\Bbb N}),E((X_n,\varphi_n)_{n\in\Bbb N})$
be two sum-like equivalence relations, both satisfying $(m1)$ and
$(m2)$. If $\psi_n,\varphi_n\le 1$, then
$$\begin{array}{cl}& E((X_n,\psi_n)_{n\in\Bbb
N})=E((X_n,\varphi_n)_{n\in\Bbb N})\cr\iff& \exists A\ge 1\forall
u,v\in X(A^{-1}\psi(u,v)\le\varphi(u,v)\le A\psi(u,v)).\end{array}$$
\end{corollary}

Before introducing the Metrization lemma, we recall several basic
notions on relations. Let $X$ be a non-empty set, we denote
$\Delta(X)=\{(u,u):u\in X\}$. A subset $U\subseteq X^2$ is called
symmetric if $(u,v)\in U\Rightarrow(v,u)\in U$ for $u,v\in X$. For
$U,V\subseteq X^2$ we define $U\circ V\subseteq X^2$ by
$$(u,r)\in U\circ V\iff\exists v((u,v)\in U,(v,r)\in V)\quad(\forall
u,v,r\in X).$$

\begin{lemma}[Metrization lemma \cite{kelley}, p. 185]
Let $U_n,\,n\in\Bbb N$ be a sequence of subsets of $X^2$ such that
\begin{enumerate}
\item[(i)] $U_0=X^2$;
\item[(ii)] each $U_n$ is symmetric and $\Delta(X)\subseteq U_n$;
\item[(iii)] $U_{n+1}\circ U_{n+1}\circ U_{n+1}\subseteq U_n$ for each $n$.
\end{enumerate}
Then there is a pseudo-metric $d$ on $X$ satisfying that
$$U_n\subseteq\{(u,v):d(u,v)<2^{-n}\}\subseteq U_{n-1}\quad(\forall
n\ge 1).$$
\end{lemma}

\begin{theorem}
Let $(X_n,\psi_n)_{n\in\Bbb N}$ satisfy $(m1)$ and $(m2)$. Then
$E((X_n,\psi_n)_{n\in\Bbb N})$ is an equivalence relation iff it is
an $\ell_p$-like equivalence relation.
\end{theorem}

\begin{proof}
Assume that $E((X_n,\psi_n)_{n\in\Bbb N})$ is an equivalence
relation. Without loss of generality, we may assume that $\psi_n\le
1$. From Lemma \ref{C}, there is $C\ge 1$ such that, for $u,v,r\in
X$,
$$\psi(v,u)\le C\psi(u,v)\mbox{ and }\psi(u,r)\le
C(\psi(u,v)+\psi(v,r)).$$ Therefore, if
$\psi(u,v),\psi(v,r),\psi(r,s)<\varepsilon$, then
$$\psi(u,s)\le C(\psi(u,v)+C(\psi(v,r)+\psi(r,s)))<(2C^2+C)\varepsilon.$$

Now denote $B=2C^2+C$. We define $U_0=X^2$ and
$$U_n=\{(u,v):\psi(u,v)<B^{-n},\psi(v,u)<B^{-n}\}\quad(n\ge 1).$$
It follows that, for each $n$, $U_n$ is symmetric and $U_{n+1}\circ
U_{n+1}\circ U_{n+1}\subseteq U_n$. By Lemma \ref{C}.(i),
$\Delta(X)\subseteq U_n$. Then the metrization lemma gives a
pseudo-metric $d$ on $X$ such that
$U_n\subseteq\{(u,v):d(u,v)<2^{-n}\}\subseteq U_{n-1}$.

It is easy to check that, $d(u,v)=0$ iff $\psi(u,v)=0$. If
$\psi(u,v)\ge B^{-2}$, then $(u,v)\notin U_2$, $d(u,v)\ge 2^{-3}$.

Denote $p=\log_2 B\ge 1$. If $0<\psi(u,v)<B^{-2}$, assume that
$B^{-(n+1)}\le\psi(u,v)<B^{-n}$ for some $n\ge 2$. Then
$\psi(v,u)<CB^{-n}<B^{-(n-1)}$. It follows that $(u,v)\in
U_{n-1},(u,v)\notin U_{n+1}$, $2^{-(n+2)}\le d(u,v)<2^{-(n-1)}$. So
$$B^{-2}d(u,v)^p<B^{-2}(2^{-(n-1)})^p\le\psi(u,v)<B^2(2^{-(n+2)})^p\le B^2d(u,v)^p.$$

Therefore, we have $E((X_n,\psi_n)_{n\in\Bbb
N})=E((X_n,d\upharpoonright X_n);p)$.
\end{proof}

\begin{corollary}
Let $(X_n,\psi_n)_{n\in\Bbb N}$ satisfy that, for $n<m$,
$X_n\subseteq X_m$ and $\psi_n=\psi_m\upharpoonright X_n^2$. Then
$E((X_n,\psi_n)_{n\in\Bbb N})$ is an equivalence relation iff it is
an $\ell_p$-like equivalence relation.

In particular, $E(X,\psi)$ is an equivalence relation iff there are
pseudo-metric $d$ on $X$ and $p\ge 1$ such that
$E(X,\psi)=E(X,d;p)$.
\end{corollary}

\section{Further remarks}

Let us consider a special case of sum-like equivalence relation. Let
$\psi_f(u,v)=f(|u-v|)$ where $f:\Bbb R^+\to\Bbb R^+$. Then
conditions (i) and (ii) for $E(\Bbb R,\psi_f)$ in Lemma \ref{C} read
as (see also \cite{matrai}, Proposition 2).
\begin{enumerate}
{\it \item[(i)] $f(0)=0$;
\item[(ii)] there is a $C\ge 1$ such that for $s,t\in\Bbb R^+$,
$$f(s+t)\le C(f(s)+f(t)),\quad f(s)\le C(f(s+t)+f(t)).$$}
\end{enumerate}

Denote $\mathcal N_f=\{x\in\Bbb R^\Bbb N:\sum_{n\in\Bbb
N}f(|x(n)|)<+\infty\}$. Then $E(\Bbb R,\psi_f)$ is an equivalence
relation iff $\mathcal N_f$ is a subgroup of $(\Bbb R^\Bbb N,+)$.
Furthermore, another interesting problem is to determine when
$\mathcal N_f$ is a linear subspace of $\Bbb R^\Bbb N$. This problem
was studied by S. Mazur and W. Orlicz (see \cite{MO}, 1.7). It was
also considered in \cite{MO} that when $\mathcal N_f$'s are Banach
spaces.

\begin{theorem}[Mazur-Orlicz] \label{MO}
Let $f:\Bbb R^+\to\Bbb R^+$ satisfy that, as $n\to\infty$,
$f(t_n)\to 0$ iff $t_n\to 0$. The necessary and sufficient condition
for $\mathcal N_f$ to be a linear space is that
\begin{enumerate}
\item[(a)] there exist constants $C>0,\varepsilon>0$ such that $f(s+t)\le C(f(s)+f(t))$
for $s,t<\varepsilon$;
\item[(b)] for every $\rho>0$ there are constants $D>0,\delta>0$
such that $f(s)\le Df(t)$ for $t<\delta,s<\rho t$.
\end{enumerate}
\end{theorem}

Note that for $t_n\ge\min\{\varepsilon,\delta\}>0,\,n\in\Bbb N$, we
have $f(t_n)\not\to 0$. Thus there is $c>0$ such that $f(t)\ge c$
for $t\ge\min\{\varepsilon,\delta\}$. If we assume that $f\le 1$,
then conditions (a) and (b) in this theorem turn to
\begin{enumerate}
{\it \item[(a)'] there exists a constant $C'>0$ such that $f(2s)\le
C'f(s)$ for $s\in\Bbb R^+$;
\item[(b)'] there exists a constant $D'>0$ such that $f(s)\le D'f(t)$ for
$s<t$.}
\end{enumerate}

For almost all known examples of sum-like equivalence relations
$E(\Bbb R,\psi_f)$, $\mathcal N_f$'s are linear spaces. In the end,
we present an example in which $\mathcal N_f$ is not linear as
follows.

\begin{example}
Let $g:\Bbb R^+\to\Bbb R^+$ be an increasing function with $g(0)=0$,
such that $g'(t)$ is decreasing with $\lim_{t\to 0}g'(t)=+\infty$.
For example, $g(x)=\sqrt{x}$ is such a function. Let
$(a_n)_{n\in\Bbb N}$ be a strictly decreasing sequence of positive
numbers with $\lim_{n\to\infty}a_n=0$.

Denote $k_n=\frac{g(a_n)}{a_n}$. Then $k_n<k_{n+1}$. We consider
equations $y=k_nx$ and $y-g(a_{n+1})=-k_{n+1}(x-a_{n+1})$. Their
solution is $(b_n,k_nb_n)$ where
$b_n=\frac{2k_{n+1}a_{n+1}}{k_n+k_{n+1}}$. We can see that
$a_{n+1}<b_n<a_n$.

Now define $f:\Bbb R^+\to\Bbb R^+$ by
$$f(t)=\left\{\begin{array}{ll}0,& t=0,\cr -k_{n+1}(t-a_{n+1})+g(a_{n+1}),&
a_{n+1}\le t<b_n,\cr k_nt,& b_n\le t<a_n,\cr g(a_0),& a_0\le
t.\end{array}\right.$$ It is easy to check that, for $s,t\in\Bbb
R^+$,
$$f(s+t)\le f(s)+f(t),\quad f(s)\le f(s+t)+f(t).$$
Note that
$\frac{f(a_{n+1})}{f(b_n)}=\frac{1}{2}\left(1+\frac{k_{n+1}}{k_n}\right)$.
Since $\lim_{t\to 0}\frac{g(t)}{t}=\lim_{t\to 0}g'(t)=+\infty$, we
can find a sequence $(a_n)_{n\in\Bbb N}$ such that
$\frac{k_{n+1}}{k_n}\to +\infty$, $\frac{f(a_{n+1})}{f(b_n)}\to
+\infty$ as $n\to\infty$. Then condition (b) in Theorem \ref{MO}
fails.
\end{example}

\end{document}